\documentclass{article}
\usepackage{emsnewsletter}

\usepackage[utf8]{inputenc} 
\usepackage[T1]{fontenc}    
\usepackage{hyperref}       
\usepackage{url}            
\usepackage{booktabs}       
\usepackage{amsfonts}       
\usepackage{nicefrac}       
\usepackage{microtype}      
\usepackage{xcolor}         
\usepackage{amsmath}
\usepackage{amssymb}
\usepackage{mathtools}
\usepackage{amsthm}
\usepackage{makecell}
\usepackage[font=small]{caption}
\usepackage{multirow}
\usepackage{makecell}
\usepackage{algorithm}
\usepackage{algorithmic}
\usepackage{lscape}
\usepackage{wrapfig}

\newcommand{\algname}[1]{#1}

\allowdisplaybreaks
\usepackage{colortbl}
\definecolor{bgcolor}{rgb}{0.8,1,1}
\definecolor{bgcolor2}{rgb}{0.8,1,0.8}
\usepackage{threeparttable}

\usepackage[capitalize,noabbrev]{cleveref}

\theoremstyle{plain}
\newtheorem{assumption}{Assumption}

\usepackage[]{todonotes}

\def\R{\mathbb{R}}

\def\C{\mathcal C}
\def\X{\mathcal X}
\def\Y{\mathcal Y}

\def\R{\mathbb R}

\def\EE{\mathbb E}
\def\PP{\mathbb P}

\def\Z{\mathcal Z}

\def\la{\langle}
\def\ra{\rangle}

\newcommand{\mA}{\mathbf{A}}

\newcommand{\eqdef}{\vcentcolon=}

\newcommand{\cO}{\mathcal{O}}

\def\<#1,#2>{\langle #1,#2\rangle}

\author{Aleksandr Beznosikov (MIPT, Moscow, Russia), Boris Polyak (ICS RAS, Moscow, Russia), Eduard Gorbunov (MIPT, Moscow, Russia), Dmitry Kovalev (KAUST, Thuwal, Saudi Arabia), Alexander Gasnikov (MIPT, IITP RAS, Moscow, Russia; Caucasus Mathematical Center, Adyghe State University, Maikop, Russia)}
\title{Smooth Monotone Stochastic Variational Inequalities and Saddle Point Problems: A Survey}

\begin{document}
   
\maketitle

\begin{abstract}
This paper is a survey of methods for solving smooth (strongly) monotone stochastic variational inequalities. To begin with, we give the deterministic foundation from which the stochastic methods eventually evolved. Then we review methods for the general stochastic formulation, and look at the finite sum setup. The last parts of the paper are devoted to various recent (not necessarily stochastic) advances in algorithms for variational inequalities. 
\end{abstract}

\section{Introduction}

In its long, more than half-century history of study \cite{stampacchia1964formes}, variational inequalities have become one of the most popular and universal optimization formulations. 
Variational inequalities are used in various areas of applied mathematics. Here we can highlight both classic examples from game theory, economics, operator theory, convex analysis \cite{arrow1958studies,browder1966existence,rockafellar1969convex,sibony1970methodes,stampacchia1964formes}, as well as newer and even more recent applications in optimization and machine learning: non-smooth optimization \cite{nesterov2005smooth}, unsupervised learning \cite{bach2008convex,esser2010general,chambolle2011first}, robust/adversarial optimization \cite{BenTal2009:book},  GANs \cite{goodfellow2014generative} and reinforcement learning \cite{Omidshafiei2017:rl,Jin2020:mdp}. 
Modern times present a new challenges to the community. The increase in scale of problems and the desire to speed up solution processes have led to a huge interest in \emph{stochastic} formulations of applied tasks, including variational inequalities. A survey of stochastic methods for solving variational inequalities is the subject of this paper.

\hspace{-0.45cm} \textbf{Structure of the paper.} In Section \ref{sec:setting}, we give a formal statement of the variational inequality problem, basic examples, and main assumptions. Section \ref{sec:determ} deals with deterministic methods, from which stochastic methods have been developed. Section \ref{sec:stoch} covers the stochastic methods. Section \ref{adv} is devoted to the recent advances in (not necessarily stochastic) variational inequalities and saddle point problems.

\section{Problem: setting and assumptions} \label{sec:setting}

\textbf{Notation.} We use $\la x,y \ra \eqdef \sum_{i=1}^nx_i y_i$ to denote standard inner product of $x,y\in\R^d$ where $x_i$ corresponds to the $i$-th component of $x$ in the standard basis in $\R^d$. It induces $\ell_2$-norm in $\R^d$ in the following way $\|x\|_2 \eqdef \sqrt{\la x, x \ra}$. We denote $\ell_p$-norms as $\|x\|_p \eqdef \left(\sum_{i=1}^d|x_i|^p\right)^{\nicefrac{1}{p}}$ for $p\in [1,\infty)$ and for $p = \infty$ we use $\|x\|_\infty \eqdef \max_{1\le i\le d}|x_i|$. The dual norm $\|\cdot\|_*$ for the norm $\|\cdot\|$ is defined in the following way: $\|y\|_* \eqdef \max\left\{\la x, y \ra\mid \|x\| \le 1\right\}$. Operator $\EE\left[\cdot\right]$ denotes full mathematical expectation. Finally, we need to introduce $\mathcal{O},\Omega$-notation to hide numerical constants which do not depend on any problem parameter, and notation $ \mathcal{\tilde O},\tilde \Omega$ to hide numerical constants and logarithmic factors.
\vspace{0.3cm}

We study variational inequalities (VI) of the  form
\begin{equation}\begin{aligned}
    \label{eq:VI}
    \text{Find} ~~ z^* &\in \Z ~~ \text{such that} ~~
    \langle F(z^*), z - z^* \rangle \geq 0, ~~ \forall z \in \Z,
\end{aligned}\end{equation}
where $F: \Z \to \R^d $ is an operator, and $\Z \subseteq \R^d$ is a convex set.

To emphasise the extensiveness of the formalism \eqref{eq:VI}, we give some examples of variational inequalities arising in applied science.

\begin{example}[Minimization] \label{ex:min} Consider the minimization problem:
\begin{align}
\label{eq:min}
\min_{z \in \Z} f(z).
\end{align}
Suppose that $F(z) \eqdef \nabla f(z)$. Then, if $f$ is convex, it can be proved that $z^* \in \Z$ is a solution for \eqref{eq:VI} if and only if $z^* \in \Z$ is a solution for \eqref{eq:min}. 
\end{example}

\begin{example}[Saddle point problem] \label{ex:minmax} Consider the saddle point problem (SPP):
\begin{align}
\label{eq:minmax}
\min_{x \in \X} \max_{y \in \Y} g(x,y).
\end{align}
Suppose that $F(z) \eqdef F(x,y) = [\nabla_x g(x,y), -\nabla_y g(x,y)]$ and $\Z = \X \times \Y$ with $\X \subseteq \R^{d_x}$, $\Y \subseteq \R^{d_y}$. Then, if $g$ is convex-concave, it can be proved that $z^* \in \Z$ is a solution for \eqref{eq:VI} if and only if $z^* \in \Z$ is a solution for \eqref{eq:minmax}. 
\end{example}

The study of saddle point problems is often associated with variational inequalities. 

\begin{example}[Fixed point problem]
Consider the fixed point problem:
\begin{align}
\label{eq:fixedpoint}
 \text{Find} ~~ z^* &\in \R^d ~~ \text{such that} ~~
    T(z^*) = z^*,
\end{align}
where $T: \R^d \to \R^d $ is an operator. With $F(z) = z - T(z)$, it can be proved that $z^* \in \Z = \R^d$ is a solution for \eqref{eq:VI} if and only if $F(z^*) = 0$, i.e. $z^* \in \R^d$ is a solution for \eqref{eq:fixedpoint}.
\end{example}

For the operator $F$ from \eqref{eq:VI} we assume the following.

\begin{assumption}[Lipschitzness] \label{as:Lipsh}
The operator $F$ is $L$-Lipschitz continuous, i.e. for all $u, v \in \Z$ we have
$
\| F(u)-F(v) \|_*  \leq L\|u - v\|.
$
\end{assumption}
In the context of \eqref{eq:min} and \eqref{eq:minmax},  $L$-Lipschitzness of the operator means that the functions $f(z)$ and $g(x,y)$ are $L$-smooth.

\begin{assumption}[Strong monotonicity]\label{as:strmon}
The operator $F$ is $\mu$-strongly monotone, i.e., for all $u, v \in \Z$ we have
$
\langle F(u) - F(v), u - v \rangle \geq \mu \| u-v\|^2_2.
$
If $\mu = 0$, then the operator $F$ is monotone.
\end{assumption}
In the context of \eqref{eq:min} and \eqref{eq:minmax}, strong monotonicity of $F$ means strong convexity of $f(z)$ and strong convexity-strong concavity of $g(x,y)$. In this paper we first concentrate on the strongly monotone and monotone cases. But there are also various assumptions relaxing monotonicity and strong monotonicity (e.g., see \cite{hsieh2020explore} and references therein).

One can point out that Assumptions \ref{as:Lipsh} and \ref{as:strmon} are sufficient for the existence of a solution \eqref{eq:VI} \cite{VIbook2003}.

Since we work on the set $\Z$, let us introduce the Euclidean projection on this set:
$$
P_{\Z}(z)=\arg\min_{v\in \Z} \|z-v\|_2.
$$

To characterise the convergence of the methods for monotone variational inequalities we introduce the gap function:
\begin{equation}
    \label{eq:gap_VI}
    \text{Gap}_{\text{VI}} (z) \eqdef \sup_{u \in \mathcal{\Z}} \left[ \langle F(u),  z - u  \rangle \right].
\end{equation}
Such a gap function, as a convergence criterion, is more suitable for the following variational inequality problem: $\langle F(z), z^* - z \rangle \leq 0$ for $z \in \Z$. Such a solution is also called weak or Minty (the solution of \eqref{eq:VI} is called strong or Stampacchia). However, according to Assumption \ref{as:Lipsh} we have that $F$ is single-valued and continuous on $\Z$, meaning that both formulations of the variational inequality problem are equivalent \cite{VIbook2003}.

For the minimization problem \eqref{eq:min}, the functional distance to the solution: $f(z) - f(z^*)$, can be used instead of \eqref{eq:gap_VI}. For saddle point problems \eqref{eq:minmax}, the gap function is also used, but it is slightly different:
\begin{equation}
    \label{eq:gap_minmax}
    \text{Gap}_{\text{SPP}} (z) \eqdef \text{gap}(x,y) = \max_{y' \in \mathcal{Y}} f(x, y') - \min_{x' \in \mathcal{X}} f(x', y).
\end{equation}
For both functions \eqref{eq:gap_VI} and \eqref{eq:gap_minmax}, it is crucial that the feasible set is bounded (in fact it is not necessary to take the whole set $\Z$ which can be unbounded, it is enough to take a bounded convex subset $\C$ which contains some solution -- see \cite{nesterov2007dual}). Therefore it is necessary to define a distance on the set $\Z$. Since this survey covers methods not only in the Euclidean setup, let us introduce a more general notion of distance.

\begin{definition}[Bregman divergence]\label{def:Bregman}
    Assume that function $\nu(z)$ is $1$-strongly convex w.r.t. $\|\cdot\|$-norm and differentiable on $\Z$ function. Then for any two points $z, z'\in \Z$ we define Bregman divergence $V(z,z')$ associated with $\nu(z)$ as follows:
    \begin{equation*}
        V(z,z') \eqdef \nu(z') - \nu(z) - \la\nabla \nu(z), z'-z \ra.
    \end{equation*}
\end{definition}
We denote the Bregman-diameter of the set $\Z$ w.r.t.\ $V(z,z')$ as $D_{\Z,V} \eqdef \max\{\sqrt{2V(z,z')}\mid z,z' \in \Z\}$. In the Euclidean case, we use $D_{\Z}$ instead of $D_{\Z,V}$. Using the definition of $V$, we denote the proximal operator as follows:
\begin{equation*}
    \text{prox}_{x}(y) = \arg\min_{z \in \mathcal{Z}}\left\{ \langle y, z \rangle + V(z,x)\right\}.
\end{equation*}

\section{Deterministic foundation: Extragradient and others} \label{sec:determ}

The first and the simplest method for solving variational inequality \eqref{eq:VI} is iterative scheme (also known as \algname{Gradient method})
\begin{equation}\label{eq:iter}
 z^{k+1} = P_{\Z} (z^k-\gamma F(z^k)),  
\end{equation}
where $\gamma>0$ is a step size.
Note that this method can be rewritten using the proximal operator with the Euclidean Bregman divergence:
\begin{equation*}
z^{k+1} = \text{prox}_{z^k}(\gamma F(z^k)).
\end{equation*}
The basic result is the convergence of the method to the unique solution of
\eqref{eq:VI} for strongly monotone and $L$-Lipschitz operator $F$; it was obtained in the papers \cite{browder1966existence,rockafellar1969convex,sibony1970methodes}.

\begin{theorem} \label{th:gm}
If Assumptions \ref{as:Lipsh}, \ref{as:strmon} hold and $0<\gamma< 2\mu/L^2$, then after $k$ iterations the method \eqref{eq:iter} converges to $z^*$ with linear rate: 
\begin{equation*}
\|z^k-z^*\|_2^2 = \mathcal{O}(R_0^2 q^k),~~\text{with}~~q = (1-2\gamma\mu+\gamma^2L^2)
\end{equation*}
and $R_0 = \|z^0 - z^*\|_2$ (here and below).
For $\gamma=\mu/L^2$ we have $q=(1-1/\kappa^2), \kappa=L/\mu$, thus the upper bound on the number of iterations to achieve the $\varepsilon$-solution (i.e., $\|z^k-z^*\|_2^2 \leq \varepsilon$) is $\mathcal{O}(\kappa^2\log(R_0^2/\varepsilon))$.
\end{theorem}

Various extensions of this statement (for $F$ being non-Lipschitz but with linear grow bounds, or for values of $F$ corrupted by noise) can be found in Theorem 1 from \cite{bakushinskii1974solution}. 

If $F$ is a potential operator (see Example \ref{ex:min}) the method \eqref{eq:iter} coincides with the gradient projection algorithm. It converges for strongly monotone $F$. Moreover, bounds for the admissible step size are less restrictive ($0<\gamma<2/L$) and complexity estimates are better ($O(\kappa\log(R_0^2/\varepsilon))$) than in Theorem \ref{th:gm}, see Theorem 2 from Section 1.4.2 of \cite{polyak1987introduction}. 

However, in the general monotone but non-strongly monotone case (for instance, for convex-concave SPP, Example \ref{ex:minmax}) convergence is lacking. Original statements on the convergence of Uzava method (a version of 
\eqref{eq:iter}) for saddle point problems
\cite{arrow1958studies} were wrong; numerous examples of divergence of the method \eqref{eq:iter} for $F$ corresponding to bilinear SPP are well known, see e.g. Figure 39 from \cite{polyak1987introduction}. 

There were many other attempts to recover convergence of gradient-like methods not for VIs, but for saddle point problems. One of them is based on the transition to modified Lagrangians when $g(x,y)$ is Lagrange function, see \cite{golshtein1977convergence,polyak1987introduction}. However, we focus on general VI case.
A possible approach is based on 
\textit{regularization} idea. Instead of the monotone variational inequality \eqref{eq:VI} one can deal with the regularized one, when monotone $F$ is replaced with strongly monotone $F+\varepsilon_k T$, where $T(z)$ is a strongly monotone operator and $\varepsilon_k>0$ is a regularization parameter. If we denote $z^k$ as the solution of regularized VI, then it is possible to prove \cite{bakushinskii1974solution} that $z^k$ converges to $z^*$ for $\varepsilon_k\rightarrow 0$. However, the solution $z^k$ usually is not easily available. To adress this problem, the  \textit{iterative regularization} technique is proposed in \cite{bakushinskii1974solution}, where one step of the basic method \eqref{eq:iter} is applied for the regularized problem. Step sizes and regularization parameters can be adjusted to guarantee the convergence.

Another technique is based on the \algname{Proximal Point}  method proposed independently by B.~Martinet in \cite{martinet1970regularisation} and by T.~Rockafellar in \cite{rockafellar1976monotone}. At each iteration it requires the solution of VI with $F+cI$, where $c>0$ and $I$ is the identity operator. This is implicit method (similar with regularization method), however there exist numerous implementable versions of \algname{Proximal Point}. For instance, some methods discussed further can be considered from this point of view.

The breakthrough for solving (non-strongly) monotone variational inequalities was made by Galina Korpelevich \cite{korpelevich1976extragradient}. She exploited the idea of the extrapolation for the gradient method. It can be explained for the simplest example of a two-dimensional min-max problem with $g(x,y)=xy, \Z=R^2$. It has the unique saddle point $z=0$, and in any point $z^k$ the direction $F(z^k)$ is orthogonal to $z^k$; thus the iteration \eqref{eq:iter} enlarges the distance to the saddle point. However, if we make the step \eqref{eq:iter} and get extrapolated point $z^{k+1/2}$, the direction $-F(z^{k+1/2})$ is attracting to the saddle point.
Thus, the \algname{Extragradient} method for solving \eqref{eq:VI} reads:
\begin{align}\label{eq:extra}
\begin{split}
    z^{k+1/2} = P_{\Z}(z^k-\gamma F(z^k)),\\
z^{k+1} = P_{\Z}(z^{k}-\gamma F(z^{k+1/2})).
\end{split}
\end{align}\label{extra}
\begin{theorem}
Let $F$ satisfy Assumptions \ref{as:Lipsh}, \ref{as:strmon} (with $\mu = 0$) and $0<\gamma<1/L$, then the Extragradient method generates a sequence of iterates $z^k$ that converge to $z^\star$.
\end{theorem}

For particular cases of the zero-sum matrix game or the general bilinear problem $g(x,y)=y^\top \mA x -b^\top x + c^{\top}y$  the method converges linearly, provided that the optimal solution is unique (see Theorem 3 from \cite{korpelevich1976extragradient}). In this case the convergence rate is equal to $\mathcal{O}(\kappa\log(R_0^2/\varepsilon))$ with $\kappa=\lambda_{\max}(\mA\mA^\top)/\lambda_{\min}(\mA \mA^\top)$. More general upper bounds for the \algname{Extragradient} method can be found in \cite{tseng1995linear} and in the recent paper \cite{mokhtari2020unified}. 
In particular, for the strongly monotone case the estimate  $O(\kappa\log(R_0^2/\varepsilon))$ with $\kappa = L/\mu$ holds true (compare with much worse bound $O(\kappa^2\log(R_0^2/\varepsilon))$ for \algname{Gradient method}). Adaptive version of the \algname{Extragradient} method (no knowledge of $L$ required) is proposed in \cite{khobotov1987modification}.

Another version of the \algname{Extragradient} method for finding saddle points is provided in \cite{korpelevich1983extrapolation}. Considering the setup of Example \ref{ex:minmax}, we can exploit just one extrapolating step for variables $y$:
\begin{align}
\label{eq:extr_spp}
y^{k+1/2}=P_{\Y}(y^k + \gamma \nabla_y g(x^k, y^k)), \nonumber\\
x^{k+1}=P_{\X}(x^k- \gamma \nabla_x g (x^k, y^{k+1/2}),\\
y^{k+1}=y^k+q (y^{k+1/2} -y^k), \nonumber
\end{align}
with $0<\gamma<1/2L$ and $0<q<1$. This method converges to the solution and if $g(x,y)$ is linear in $y$, then the convergence rate is linear. If we put $q = 1$ in the method \eqref{eq:extr_spp}, then $y^{k+1} = y^{k+1/2}$ and we get the so-called Alternating Gradient method (Alternating descent-ascent). In \cite{zhang2022near}, it was proved  that this method has \textit{local} linear convergence with complexity $O(\kappa\log(R_0^2/\varepsilon))$, where $\kappa = L/\mu$.

L.~Popov \cite{popov1980modification} proposed a version of extrapolation scheme (sometimes this type of scheme is called \textit{optimistic} or \textit{single-call}):
\begin{align}\label{eq:popov}
\begin{split}
    z^{k+1/2} = P_{\Z}(z^k-\gamma F(z^{k-1/2})),\\
z^{k+1} = P_{\Z}(z^{k}-\gamma F(z^{k+1/2})).
\end{split}
\end{align}
It requires the single calculation of $F$ at each iteration vs two calculations in the \algname{Extragradient} method. As shown in
\cite{popov1980modification}, the method \eqref{eq:popov} converges for $0<\gamma<1/3L$. Convergence rates for this method were obtained recently \cite{gidel2018variational, mokhtari2020unified}, it is $O(\kappa\log(R_0^2/\varepsilon))$ with $\kappa = L/\mu$ for the strongly monotone case and $\kappa=\lambda_{\max}(\mA\mA^\top)/\lambda_{\min}(\mA\mA^\top)$ for the bilinear case. It is important to note that in the general strongly monotone case this estimate is optimal \cite{zhang2021lower}, but for the bilinear problem the upper bounds available in the literature for both the \algname{Extragradient} and optimistic methods are not tight \cite{ibrahim2020linear}. Meanwhile, optimal estimates $O(\sqrt{\kappa}\log(R_0^2/\varepsilon))$ with $\kappa = \lambda_{\max}(\mA\mA^\top)/\lambda_{\min}(\mA\mA^\top)$ can be achieved using approaches from \cite{azizian2020accelerating} and \cite{alkousa2020accelerated}.

The extension of the above schemes to an arbitrary proximal setup was obtained in the work of A.~Nemirovsky \cite{nemirovski2004prox}. He proposed the \algname{Mirror-Prox} method for VIs, exploiting Bregman divergence:
\begin{align}\label{eq:MirrorProx}
    \begin{split}
        z^{k+1/2} &= \text{prox}_{z^{k}}\left(\gamma F(z^k)\right), \\
    z^{k+1} &= \text{prox}_{z^{k}}\left(\gamma F(z^{k+1/2})\right).
    \end{split}
\end{align}
It implies the following result on the convergence rate.
\begin{theorem}\label{th:MirrProx}
Let $F$ satisfy Assumptions \ref{as:Lipsh}, \ref{as:strmon} (with $\mu = 0$) and 
\begin{equation}\label{nemir}
    \hat{z}^k=\frac{1}{k}\sum_{i=1}^k z^{i+1/2},
\end{equation}
where $z^{i+1/2}$ are generated by algorithm \eqref{eq:MirrorProx} with $\gamma=1/\sqrt{2}L$. Then, after $k$ iterations
\begin{equation}\label{eq:MP}
\text {Gap}_{VI} ( \hat{z}^k) = \cO \left( \frac{LD_{\Z,V}^2}{k} \right).
\end{equation}
\end{theorem}
Numerous extensions of these original versions of iterative methods for solving variational inequalities were published later. One can highlight Tseng's \algname{Forward-Backward Splitting} \cite{tseng2000modified}, Nesterov's \algname{Dual Extrapolation} \cite{nesterov2007dual}, Malitsky and Tam's \algname{Forward-Reflected-Backward} \cite{doi:10.1137/18M1207260}. All methods have convergence guarantees \eqref{eq:MP}. It turns out that this rate is optimal \cite{ouyang2021lower}.


\section{Stochastic methods: different setups and assumptions} \label{sec:stoch}

In this section, we move from deterministic to stochastic methods, i.e., we consider \eqref{eq:VI} with the following operator:
\begin{equation}
    \label{eq:stoch_VI}
    F(z) = \EE_{\xi \sim \mathcal{D}} \left[F_{\xi}(z)\right],
\end{equation}
where $\xi$ is a random variable, $\mathcal{D}$ is some (typically unknown) distribution, $F_{\xi}: \Z \to \R^d $ is a stochastic operator. In this setup, calculating the value of the full operator $F$ is computationally expensive or intractable. Therefore, one has to work mainly with stochastic realizations $F_{\xi}$.

\subsection{General case}\label{sec:general_case}

The stochastic formulation \eqref{eq:stoch_VI} for the problem \eqref{eq:VI} was first considered by the authors of \cite{juditsky2011solving}. They proposed a natural stochastic generalization of \algname{Extragradient} (more precisely, of \algname{Mirror-Prox}):
\begin{align}
\label{eq:iseg}
    \begin{split}
    z^{k+1/2} &= \text{prox}_{z^{k}}\left(\gamma F_{\xi^k} (z^k)\right), \\
    z^{k+1} &= \text{prox}_{z^{k}}\left(\gamma F_{\xi^{k+1/2}} (z^{k+1/2})\right).
    \end{split}
\end{align}
Here it is important to note that $\xi^k$ and $\xi^{k+1/2}$ are independent and $F_{\xi}(z)$ is an unbiased estimator of $F(z)$. Moreover, $F_{\xi}(z)$ is assumed to satisfy the following condition.
\begin{assumption}[Bounded variance]\label{as:bounded_variance}
The unbiased operator $F_{\xi}$ has uniformly bounded variance, i.e., for all $\xi \sim \mathcal{D}$ and $u \in \Z$ we have
$
\EE\|F_{\xi}(u) - F(u) \|^2_* \leq \sigma^2.
$
\end{assumption}
Under this assumption, the next result was derived in \cite{juditsky2011solving}.
\begin{theorem} \label{th:iseg_mon}
Let $F_{\xi}$ satisfy Assumptions \ref{as:Lipsh}, \ref{as:strmon} (with $\mu = 0$), \ref{as:bounded_variance} and $\hat{z}^k$ be defined as in \eqref{nemir}, where $z^{i+1/2}$ are generated by the algorithm \eqref{eq:iseg} with $\gamma= \min\left\{ \frac{1}{\sqrt{3}L}, D_{\Z,V} \sqrt{\frac{1}{7k\sigma^2}}\right\}$. Then, after $k$ iterations we can guarantee that 
\begin{align}
    \label{eq:iseq_mon}
    \EE\left[ \text {Gap}_{VI} ( \hat{z}^k) \right] = \cO\left( \frac{L D^2_{\Z,V}}{k} +  D_{\Z,V}\sqrt{\frac{\sigma^2}{k}} \right).
\end{align}
\end{theorem}
In \cite{beznosikov2020distributed}, the authors gave an analysis of the algorithm \eqref{eq:iseg} for strongly monotone VIs in the Euclidean case. In particular, under Assumptions \ref{as:Lipsh}, \ref{as:strmon}, \ref{as:bounded_variance} one can guarantee that after $k$ iterations of the method \eqref{eq:iseg} it holds (here and below we omit numerical constants in the exponent multiplier)
\begin{align}
    \label{eq:iseg_str_mon}
    \EE \|z^k - z^*\|^2_2 = \mathcal{\tilde O} \left( R_0^2 \exp \left(- \frac{\mu k }{L} \right) + \frac{\sigma^2}{\mu^2 k} \right).
\end{align}
Also, in \cite{beznosikov2020distributed}, the authors obtained lower complexity bounds for solving VIs satisfying Assumptions \ref{as:Lipsh}, \ref{as:strmon}, \ref{as:bounded_variance} with stochastic methods. 
It turns out that the results of Theorem \ref{th:iseg_mon} in the monotone case and those from \eqref{eq:iseg_str_mon} are optimal and meet lower bounds up to numerical constants.

Optimistic-like (or single-call) methods were also investigated in the stochastic setting. The work \cite{gidel2018variational} applies the following update scheme:
\begin{align}
\label{eq:past_seg}
    \begin{split}
    z^{k+1/2} &= P_{\Z}\left(z^k - \gamma F_{\xi^{k-1/2}} (z^{k-1/2})\right), \\
    z^{k+1} &= P_{\Z}\left(z^k - \gamma F_{\xi^{k+1/2}} (z^{k+1/2})\right).
    \end{split}
\end{align}
For this method, in the monotone Euclidean case, the authors proved an estimate similar to \eqref{eq:iseq_mon}. In the strongly monotone case, \eqref{eq:past_seg} was investigated in the paper \cite{hsieh2019convergence}, but their estimates do not meet the lower bounds. The optimal estimates for this scheme were obtained later in the work \cite{beznosikov2022unified}.

The work \cite{kotsalis2020simple1} deals with a slightly different single-call approach in the non-Euclidean case:
\begin{align}
\label{eq:lan_seg}
    \begin{split}
    z^{k+1} &= \text{prox}_{z^k}\left(\gamma_k F_{\xi^{k}} (z^{k}) + \gamma_k \alpha_k \left[F_{\xi^{k}} (z^{k}) - F_{\xi^{k-1}} (z^{k-1})\right] \right).
    \end{split}
\end{align}
This update is a modification of the \algname{Forward-Reflected-Backward} approach, namely here $\alpha_k$ is a parameter, while in \cite{doi:10.1137/18M1207260}, $\alpha_k \equiv 1$. The analysis of the method \eqref{eq:lan_seg} gives optimal estimates in both the strongly monotone and monotone cases.

The theoretical results and guarantees discussed above significantly rely on the bounded variance assumption (Assumption~\ref{as:bounded_variance}). This assumption is quite restrictive (especially when the domain is unbounded) and does not hold for many popular machine learning problems. Moreover, one can even design a strongly monotone variational inequality on an unbounded domain such that the method \eqref{eq:iseg} \emph{diverges} exponentially fast \cite{chavdarova2019reducing}.
Authors of \cite{hsieh2020explore, gorbunov2022stochastic} consider a relaxation of the bounded variance condition and assume that $
\EE\|F_{\xi}(u) - F(u) \|^2_2 \leq \sigma^2 + \delta \|u - z^* \|_2^2$ with $\delta \geq 0$ in the Euclidean case.
Under this condition and Assumptions \ref{as:Lipsh}, \ref{as:strmon}, the authors of \cite{gorbunov2022stochastic} proved that after $k$ iterations of the algorithm \eqref{eq:iseg} (when $\Z = \R^d$) it holds that
\begin{align}
    \label{eq:seg_noise1}
    \EE \|z^k - z^*\|^2_2 = \mathcal{O} \left( \kappa R^2_0 \exp \left(- \frac{k}{\kappa} \right) + \frac{\sigma^2}{\mu^2 k} \right),
\end{align}
where $\kappa = \max\left\{ \frac{\delta}{\mu^2}; \frac{L + \sqrt{\delta}}{\mu}\right\}$. The same assumption on stochastic realizations was considered in \cite{kotsalis2020simple}. The authors used the method \eqref{eq:lan_seg} and provided the following estimate:
\begin{align}
\label{eq:seg_noise2}
    \EE\|z^k - z^*\|^2_2 = \mathcal{O} \left( R^2_0 \exp \left(- \frac{\mu k}{L} \right) + \frac{\sigma^2 + \delta^2 D^2_{\Z}}{\mu^2 k} \right),
\end{align}
Results \eqref{eq:seg_noise1} and \eqref{eq:seg_noise2} are competitive. The first estimate is better in terms of the stochastic term (second term), while the second result is more competitive in terms of the deterministic term (first term). However, both of these results do not fully cover the issue of bounded noise, because the condition considered above is not general. The key to avoiding the assumption about the bounded variance of $F_{\xi}$ lies in the way how stochasticity is generated in the method \eqref{eq:iseg}. The method \eqref{eq:iseg} is sometimes called Independent Samples Stochastic Extragradient (\algname{I-SEG}). To address bounded variance issue, K.~Mishchenko et al. \cite{mishchenko2020revisiting} proposed another stochastic modification of the \algname{Extragradient} algorithm called Same Sample Stochastic Extragradient (\algname{S-SEG}):
\begin{align*}
    z^{k+1/2} &= z^{k} - \gamma F_{\xi^k} (z^k), \\
    z^{k+1} &= z^{k} - \gamma F_{\xi^{k}} (z^{k+1/2}).
\end{align*}
For simplicity we present the above method for the case when $\Z = \R^d$ ($F(x^*) = 0$), while \cite{mishchenko2020revisiting} contains a more general case of regularized VIs. In contrast to \algname{I-SEG}, \algname{S-SEG} uses the same sample $\xi^k$ for both steps at iteration $k$. Although such a strategy cannot be implemented in some scenarios (streaming oracle), it can be applied to finite-sum problems that have been gaining an increasing attention in the recent years. Moreover, \algname{S-SEG} significantly relies on the following assumption.
\begin{assumption}\label{as:F_xi_Lip}
    The operator $F_{\xi}(z)$ is $L$-Lipschitz and $\mu$-strongly monotone almost surely in $\xi$, i.e., $\|F_{\xi}(z) - F_\xi(z')\|_2 \leq L\|z - z'\|_2$ and $\langle F_{\xi}(z) - F_\xi(z'), z-z' \rangle \geq \mu\|z - z'\|_2^2$ for all $z,z' \in \R^d$ almost surely in $\xi$.
\end{assumption}
The evident difference between the setups for \algname{I-SEG} and \algname{S-SEG} can be explained through the connection between \algname{Extragradient} and the \algname{Proximal Point} method (\algname{PP}) \cite{martinet1970regularisation, rockafellar1976monotone}. We assume that $\Z = \R^d$ ($F(z^*) = 0$) in all results discussed further in Section~\ref{sec:general_case}. In this setup, \algname{PP} has the following update rule
\begin{equation*}
    z^{k+1} = z^k - \gamma F(z^{k+1}).
\end{equation*}
The method converges for any monotone operator $F$ and any $\gamma > 0$. However, the update rule of \algname{PP} is implicit and cannot be computed efficiently in many situations. The \algname{Extragradient} method can be seen as a natural approximation of \algname{PP} that substitutes $z^{k+1}$ in the right-hand side by one gradient step from $z^k$:
\begin{equation*}
    z^{k+1} = z^k - \gamma F\left(z^k - \gamma F(z^k)\right).
\end{equation*}
In addition, when $F$ is $L$-Lipschitz, one can estimate how good the approximation is. Consider $z^{k+1} = z^k - \gamma F(z^k - \gamma F(z^k))$ (the \algname{Extragradient} step) and $\tilde z^{k+1} = z^k - \gamma F(\tilde z^{k+1})$ (the \algname{PP} step). Then, $\|z^{k+1} - \tilde{z}^{k+1}\|_2$ can be estimated as follows \cite{mishchenko2020revisiting}:
\begin{align*}
   & \|z^{k+1} - \tilde z^{k+1}\|_2 = \gamma \|F(z^k - \gamma F(z^k)) - F(\tilde z^{k+1})\|_2 \\
    &\leq \gamma L \|z^k - \gamma F(z^k) - \tilde z^{k+1}\|_2 = \gamma^2 L \|F(z^k) - F(\tilde z^{k+1})\|_2\\
    &\leq \gamma^2 L^2 \|z^k - \tilde z^{k+1}\|_2 = \gamma^3 L^2 \|F(\tilde z^{k+1})\|_2 \\
    &\leq \gamma^3 L^3 \|\tilde z^{k+1} - z^*\|_2.
\end{align*}
That is, the difference between the \algname{Extragradient} and \algname{PP} steps is of the order $\cO(\gamma^3)$ rather than $\cO(\gamma^2)$. Since the later corresponds to the difference between \algname{PP} and simple gradient step \eqref{eq:iter}, \algname{Extragradient} better approximates \algname{PP} than gradient steps, which are known to be non-convergent for general monotone Lipschitz variational inequalities. This approximation feature of \algname{Extragradient} is crucial for its convergence and, as the above derivation implies, the approximation argument significantly relies on the Lipschitzness of operator $F$.

Let us go back to the differences between \algname{I-SEG} and \algname{S-SEG}. In \algname{S-SEG}, iteration $k$ can be considered as a single \algname{Extragradient} step for operator $F_{\xi^k}(z)$. Therefore, Lipschitzness and monotonicity of $F_{\xi^k}(z)$ (Assumption~\ref{as:F_xi_Lip}) are important for the analysis of \algname{S-SEG}. In contrast, \algname{I-SEG} uses different operators for the extrapolation and update steps. In this case, there is no effect from the Lipschitzness/monotonicity of individual $F_{\xi}(z)$. Therefore, the analysis of \algname{I-SEG} naturally relies on the Lipschitzness and monotonicity of $F(z)$ and closeness (on average) of $F_\xi(z)$ and $F(z)$ (Assumption~\ref{as:bounded_variance}).

The convergence of \algname{I-SEG} was discussed earlier in this section. Regarding \algname{S-SEG}, the following result holds \cite{mishchenko2020revisiting}.

\begin{theorem}
    Let Assumption~\ref{as:F_xi_Lip} hold. Then, there exists a choice of step size $\gamma$ \cite{gorbunov2022stochastic} such that the output of \algname{S-SEG} after $k$ iterations satisfies
    \begin{equation*}
        \EE\|z^k - z^*\|_2^2 = \cO\left(\frac{L R^2_0}{\mu} \exp\left(- \frac{\mu k}{L}\right) + \frac{\sigma_{*}^2}{\mu^2 k}\right),
    \end{equation*}
    where $\sigma_*^2 = \EE\|F_\xi(z^*)\|_2^2$.
\end{theorem}

The rate is similar to the one known for \algname{I-SEG} up to the following differences. First, instead of the uniform bound on the variance $\sigma^2$, the rate depends on $\sigma_*^2$, which is the variance of $F_\xi$ measured at the solution. In many cases, $\sigma^2 = \infty$ while $\sigma_*^2$ is finite. From this perspective, \algname{S-SEG} has better rate than \algname{I-SEG}. However, it comes with a price: while the rate of \algname{I-SEG} depends on the Lipschitz and strong-monotonicity constants of $F$, the rate of \algname{S-SEG} depends on \emph{the worst} constants of $F_\xi$ that can be much worse than those for $F$. In particular, consider the finite-sum setup with uniform sampling of $\xi$: $F(x) = \tfrac{1}{n}\sum_{i=1}^n F_i(x)$ where $F_i$ is $L_i$-Lipschitz and $\mu_i$-strongly monotone, $\PP\{\xi = i\} = \nicefrac{1}{n}$. Then, Assumption~\ref{as:F_xi_Lip} holds with $L = \max_{i\in [n]}L_i$, $\mu = \min_{i \in [n]}\mu_i$ and these constants appear in the rate from Theorem~\ref{as:bounded_variance}. The authors of \cite{gorbunov2022stochastic} tighten this rate. In particular, they prove that for \algname{S-SEG} with different step sizes for extrapolation and update steps it holds
\begin{equation*}
    \EE\|z^k - z^*\|_2^2 = \cO\left(\frac{L R^2_0}{\overline{\mu}}\exp\left(- \frac{\overline{\mu} k}{L}\right) + \frac{\sigma_{*}^2}{\overline{\mu}^2 k}\right),
\end{equation*}
where $\sigma_*^2 = \frac{1}{n}\sum_{i=1}^n\|F_i(z^*)\|_2^2$ and $\overline{\mu} = \tfrac{1}{n}\sum_{i=1}^n \mu_i$. Since $\overline{\mu}$ is (sometimes significantly) larger than $\mu$, the improvement is noticeable. Moreover, when $\{L_i\}_{i=1}^n$ are known, one can consider so-called \emph{importance sampling} \cite{gower2019sgd}: $\PP\{\xi = i\} = \nicefrac{L_i}{(n\overline{L})}$, where $\overline{L} = \tfrac{1}{n}\sum_{i=1}^n L_i$. As the authors of \cite{gorbunov2022stochastic} show, importance sampling can be combined with \algname{S-SEG} via allowing the extrapolation and update step sizes at iteration $k$ to depend on the sample $\xi^k$. In particular, for the proposed modification of \algname{S-SEG} they derive
\begin{equation*}
    \EE\|z^k - z^*\|_2^2 = \cO\left(\frac{\overline{L} R_0^2}{\overline{\mu}}\exp\left(- \frac{\overline{\mu} k}{\overline{L}}\right) + \frac{\hat\sigma_{*}^2}{\overline{\mu}^2 k}\right),
\end{equation*}
where $\hat \sigma_*^2 = \frac{1}{n}\sum_{i=1}^n\tfrac{\overline{L}}{L_i} \|F_i(z^*)\|_2^2$. The exponentially decaying term is always better than the corresponding one for \algname{S-SEG} with uniform sampling. This usually implies faster convergence during the initial stage. Next, typically, larger norm of $F_i(z^*)$ implies larger $L_i$, e.g., $\|F_i(z^*)\|_2^2 \sim L_i^2$. In this case, $\hat \sigma_*^2 \leq \sigma_*^2$, because $\hat \sigma_{*}^2 \sim (\overline{L})^2$ and $\sigma_{*}^2 \sim \overline{L^2} = \frac{1}{n}\sum_{i=1}^n L_i^2 \geq (\overline{L})^2$. Moreover, one can allow other sampling strategies and cover the case when some $\mu_i$ are negative, see \cite{gorbunov2022stochastic} for the details.

\subsection{Finite-sum case}

As noted earlier, when we deal with problem \eqref{eq:stoch_VI}, it is often the case (especially in practical problems) that the distribution $\mathcal{D}$ is unknown, but we have some samples from $\mathcal{D}$. Then, one can replace problem \eqref{eq:stoch_VI} by a finite-sum approximation:
\begin{equation}
    \label{eq:fs}
    F(z) = \frac{1}{n} \sum\limits_{i=1}^n F_i(z).
\end{equation}
This approximation is sometimes also called Monte Carlo approximation. For machine learning problems the term empirical risk is often encountered. Although calls of the full operator is now tractable, they remain expensive in practice. Therefore, it is worth to avoid frequent computing of $F$ and mainly use calls to single $F_i$ operators or small batches of them.

Before presenting the results, let us introduce the analogue of the Lipschitzness assumption.
\begin{assumption}[Lipschitzness in mean] \label{as:Lipsh_mean}
The operator $F$ is $L_{\text{avg}}$-Lipschitz continuous in
mean, i.e. for all $u, v \in \Z$ we have
$
\EE\left[\| F_{\xi}(u)-F_{\xi}(v) \|^2_*\right]  \leq L^2_{\text{avg}} \|u - v\|^2.
$
\end{assumption}
For example, if $F_i$ is $L_i$-Lipschitz for all $i$ and we draw index $\xi = i$ with probability $p_i = L_i / \sum_j L_j$, then $L_{\text{avg}} = \tfrac{1}{n} \sum_j L_j$.

The study of finite-sum problems in stochastic optimization is connected, first of all, with classical methods for minimization problems such as \algname{SVRG} \cite{NIPS2013_ac1dd209} and \algname{SAGA} \cite{NIPS2014_ede7e2b6}. For the saddle point problems, these methods were adopted in \cite{palaniappan2016stochastic} (in fact, these results are also valid for variational inequalities). The authors considered strongly convex--strongly concave saddles in the Euclidean case and proved the following estimates for \algname{SVRG} and \algname{SAGA}:
\begin{equation*}
    \EE\|z^k - z^*\|_2^2 = \cO\left( R^2_0 \exp\left(- \min \left\{\frac{1}{n}; \frac{\mu^2}{L^2_{\text{avg}}} \right\} k \right)\right).
\end{equation*}
Since the bound above is not tight in terms of $L_{\text{avg}}/\mu$, the authors proposed accelerating \algname{SVRG} and \algname{SAGA} via the Catalyst envelope \cite{lin2015universal}. In this case, they have the following bound:
\begin{equation}
    \label{eq:bach_acc}
    \EE\|z^k - z^*\|_2^2 = \cO\left( R^2_0 \exp\left(- \min \left\{\frac{1}{n}; \frac{\mu}{\sqrt{n} L_{\text{avg}}} \right\} \frac{k}{\log [L_{\text{avg}}/\mu]} \right)\right).
\end{equation}
The same estimates for saddle point problems methods based on accelerating envelopes were also presented in \cite{tominin2021accelerated}. 

An important step in the study of the finite-sum stochastic setup was the work \cite{carmon2019variance}. It is primarily focused on bilinear games. For this class of problems, the authors improved the estimate \eqref{eq:bach_acc} and removed the additional logarithmic factor. For general problems (saddle point and variational inequalities) the results of this paper are very similar to \eqref{eq:bach_acc} and also have an additional logarithmic factor. Meanwhile, the authors also considered the convex-concave/monotone case in the non-Euclidean setting and got that for their method after $k$ iterations it holds
\begin{align}
    \label{eq:carmon}
    \EE\left[ \text {Gap}_{VI} ( \hat{z}^k) \right] = \mathcal{\tilde O} \left( \frac{\sqrt{n} L_{\text{avg}} D^2_{\Z,V}}{k} \right).
\end{align}
The problem of the additional logarithmic factor was resolved in \cite{alacaoglu2021stochastic}. The authors proposed a modification of \algname{Extragradient}:
\begin{align}
\label{eq:alacouglu_seg}
    \begin{split}
    z^{k+1/2} =& P_{\Z}\left( z^k + \tau (w^k - z^k)- \gamma F(w^k) \right), \\
    \Delta^k =& F_{\xi^k}(z^{k+1/2}) - F_{\xi^k} (w^{k}) + F(w^k),\\
    z^{k+1} =& P_{\Z}\left( z^k + \tau (w^k - z^k)- \gamma \Delta^k \right), \\
    w^{k+1} =& 
    \begin{cases}
			z^{k+1},& \text{with probability }p\\
			w^{k},& \text{with probability }1-p
	\end{cases}.
    \end{split}
\end{align}
This algorithm is a combination of the extra step technique from the VIs theory and the loopless approach \cite{pmlr-v117-kovalev20a} for finite-sum problems. An interesting detail of the method is the randomized negative momentum: $\tau (w^k - z^k)$. While for minimization problems it is usual to apply positive/heavy-ball momentum, the opposite approach turns out to be useful for saddle point problems and variational inequalities. This effect was noticed earlier \cite{gidel2019negative, yoon2021accelerated, Yura2021} and appeared now in the theory of stochastic methods for VIs. Also, in \cite{alacaoglu2021stochastic}, the authors presented modifications for \algname{Forward-Backward}, \algname{Forward-Reflected-Backward} as well as for \algname{Extragradient} in the non-Euclidean case. 

As we noted earlier, the results of \cite{alacaoglu2021stochastic} give estimates \eqref{eq:bach_acc} and \eqref{eq:carmon}, but without additional logarithmic factors. That is, to achieve $\EE\|z^k - z^*\|_2^2 \leq \varepsilon$ in the strongly monotone case and $\EE[ \text {Gap}_{VI} ( \hat{z}^k) ] \leq \varepsilon$ in the monotone case the methods from \cite{alacaoglu2021stochastic} require 
\begin{equation}\label{eq:ndjkvjfdvkdnjkf}
    \cO\left(\max\left\{n; \frac{\sqrt{n}L_{\text{avg}}}{\mu}\right\}\log\frac{R_0^2}{\varepsilon}\right)
\end{equation}
and 
\begin{equation}\label{eq:vdlmdfmvo}
    \cO\left(\frac{\sqrt{n} L_{\text{avg}} D^2_{\Z,V}}{\varepsilon} \right)
\end{equation}
stochastic oracle calls respectively. It remains to discuss the effect of batching on the method from \eqref{eq:alacouglu_seg}, i.e., how the oracle complexity bounds change if we use not a single sample $F_{\xi^k}$ at each iteration, but a batch size of $b$: $\tfrac{1}{b}\sum_{i \in S^k} F_{i}$, where $S_k \subseteq [n]$ is the $b$-element set of the indices in the mini-batch. In this case, the methods from \cite{alacaoglu2021stochastic} give estimates \eqref{eq:ndjkvjfdvkdnjkf} and \eqref{eq:vdlmdfmvo}, but multiplied by an additional factor $\sqrt{b}$. This extra multiplier issue was resolved in \cite{kovalev2022optimal} using the following method:
\begin{align}
\label{eq:kovalev_seg}
    \begin{split}
    \Delta^k =& \tfrac{1}{b}\sum_{i \in S^k} \big[
    F_{i}(z^k) - F_{i} (w^{k-1}) \\
    &+ \alpha(F_{i}(z^k) - F_{i}(z^{k-1})) \big] + F(w^{k-1}), \\
    z^{k+1} =& P_{\Z}\left( z^k + \tau (w^k - z^k)- \gamma \Delta^k \right), \\
    w^{k+1} =& 
    \begin{cases}
			z^{k+1},& \text{with probability }p\\
			w^{k},& \text{with probability }1-p
	\end{cases}.
    \end{split}
\end{align}
The authors proved that in the strongly monotone case this method gives an estimate \eqref{eq:ndjkvjfdvkdnjkf}, i.e., without additional logarithmic factors and without factors depending on $b$.

The only issue that remains to be understood is whether the current state of the art methods with best complexities from \cite{alacaoglu2021stochastic, kovalev2022optimal} are optimal. The lower bounds from \cite{han2021lower} claim that under Assumptions \ref{as:Lipsh_mean} and \ref{as:strmon}, the methods above are optimal. However, under $L_{\max}$-Lipschitzness of all $F_i$, $i \in [n]$ and Assumption~\ref{as:strmon}, the lower bound from \cite{han2021lower} is
\begin{equation*}
    \EE\|z^k - z^*\|_2^2 = \Omega \left( R_0^2 \exp\left(- \min \left\{\frac{1}{n}; \frac{\mu}{ L_{\max}} \right\}k \right)\right).
\end{equation*}
The question whether this lower bound is tight remains open.

\subsection{Cocoercivity assumption}

In some papers, the following assumption is used instead of  \ref{as:Lipsh}.
\begin{assumption}[Cocoercivity] \label{as:coco}
The operator $F$ is $l$-cocoercive, i.e., for all $u, v \in \Z$ we have
$
\| F(u)-F(v) \|_2^2  \leq l \langle F(u) - F(v), u - v \rangle.
$
\end{assumption}
This assumption is stronger than monotonicity + Lipschitzness, i.e., not all monotone Lipschitz operators are cocoercive. One can note that the operator for the  bilinear SPP ($\min_x \max_y x^\top A y$) is not cocoercive. But if it is known that $F$ is $L$-Lipschitz and $\mu$-strongly monotone, then it is $L^2/\mu$-cocoercive. Moreover, if we consider a convex $L$-smooth minimization problem, then the corresponding operator is $L$-cocoercive.

There is no need to use an \algname{Extragradient} for cocoercive operators. One can apply the iterative scheme \eqref{eq:iter} and its modifications for the stochastic case. In spite of this, the first work on cocoercive operators in the stochastic cases, used the \algname{Extragradient} as the basic method \cite{chavdarova2019reducing}. In this paper, the authors investigated methods for finite-sum problems. The latter results from \cite{loizou2021stochastic, beznosikov2022stochastic} give an almost complete picture of stochastic algorithms based on method \eqref{eq:iter} for operators under Assumption \ref{as:coco}. In particular, the work \cite{beznosikov2022stochastic} gives a unified analysis for a large number of popular stochastic methods known yet for minimization problems \cite{gorbunov2020unified}.

\subsection{High-probability convergence}

Before this section, we focused on in-expectation convergence guarantees for the stochastic methods, i.e., bounds on $\EE[ \text {Gap}_{VI} ( \hat{z}^k) ]$ and/or $\EE\|z^k - z^*\|_2^2$. However, \emph{high-probability convergence guarantees}, i.e., bounds on $\text {Gap}_{VI} ( \hat{z}^k)$ and/or $\|z^k - z^*\|_2^2$ that hold with probability at least $1- \beta$ for given confidence level $\beta \in (0,1)$, reflect the real behavior of the methods more accurately \cite{gorbunov2020stochastic}. Despite this fact, high-probability convergence of stochastic methods for solving VIs is studied only in a couple of works.

It is worth mentioning that one can always deduce the high-probability bound from the in-expectation one via Markov's inequality. However, in this case, the derived rate of convergence will have a negative-power dependence on $\beta^{-1}$. Such guarantees are not desirable and the goal is to derive the rates that have (poly-)logarithmic dependence on the confidence level, i.e., $\beta$ should appear only in the $\cO(\text{poly}(\log(\nicefrac{1}{\beta})))$ factor.

The first and for many years the only high-probability guarantees of this type for solving stochastic VIs are derived in \cite{juditsky2011solving}. The authors assume that $F$ is monotone and $L$-Lipschitz, the domain is bounded, and $F_\xi$ is an unbiased estimator with sub-Gaussian (light) tails of the distribution:
\begin{equation*}
    \EE\left[\exp\left(\frac{\|F_\xi (x) - F(x)\|_2^2}{\sigma^2}\right)\right] \leq \exp(1).
\end{equation*}
The above condition is much stronger than Assumption~\ref{as:bounded_variance}. Under these assumptions, the authors of \cite{juditsky2011solving} prove that after $k$ iterations of \algname{Mirror-Prox} with probability at least $1 - \beta$ (for any $\beta \in (0,1)$) the following inequality holds:
\begin{equation*}
    \text {Gap}_{VI} ( \hat{z}^k) = \cO\left(\frac{L D_{\Z}^2}{k} + \frac{\sigma D_{\Z}\log(\nicefrac{1}{\beta})}{\sqrt{k}}\right).
\end{equation*}
Up to the logarithmic factor this result coincides with in-expectation one and, thus, it is optimal (up to the logarithms). However, the result is derived under restrictive light-tails assumption.

This limitation was recently addressed in \cite{gorbunov2022clipped}, where the authors derived the high-probability rates for the considered problem under just bounded variance assumption. In particular, they consider \algname{clipped-SEG} for problems with $\Z = \R^d$:
\begin{align*}
    z^{k+1/2} &= z^{k} - \gamma \cdot \text{clip}(F_{\xi^k} (z^k), \lambda_k), \\
    z^{k+1} &= z^{k} - \gamma \cdot \text{clip}(F_{\xi^{k+1/2}} (z^{k+1/2}), \lambda_{k+1/2}),
\end{align*}
where $\text{clip}(x,\lambda) = \min\{1, \nicefrac{\lambda}{\|x\|_2}\}x$ is the clipping operator -- a popular tool in deep learning \cite{pascanu2013difficulty, Goodfellow-et-al-2016}. In the setup, when $F$ is monotone and $L$-Lipschitz and Assumption~\ref{as:bounded_variance} holds, the authors of \cite{gorbunov2022clipped} prove that after $k$ iterations of \algname{clipped-SEG} with probability at least $1 - \beta$ (for any $\beta \in (0,1)$) the following inequality holds:
\begin{equation*}
    \text {Gap}_{VI} ( \hat{z}^k) = \cO\left(\frac{LR_0^2 \log(\nicefrac{k}{\beta})}{k} + \frac{\sigma R_0\sqrt{\log(\nicefrac{k}{\beta})}}{\sqrt{k}}\right),
\end{equation*}
Up to the differences in logarithmic factors, the definition of $\sigma$, and the difference between $D_{\Z}$ and $R_0$ the rate coincides with the one from \cite{juditsky2011solving} while derived without the light-tails assumption. The key algorithmic tool that allows removing the light-tails assumption is clipping: with a proper choice of the clipping level $\lambda$ the authors cut heavy tails without making the bias too large. It is worth mentioning that the result for \algname{clipped-SEG} is derived for the unconstrained case and the rate depends on $R_0$, while in \cite{juditsky2011solving}, the analysis relies on the boundedness of the domain, which diameter explicitly appears in their rate. To remove the dependence on the diameter of the domain, the authors of \cite{gorbunov2022stochastic} show that with the high-probability the iterates produced by \algname{clipped-SEG} stay in the ball around $x^*$ with the radius proportional to $R_0$. Using this trick, the authors of \cite{gorbunov2022stochastic} also show that it is sufficient to assume everything (monotonicity and Lipschitzness of $F$ and bounded variance) only on this ball. Such a generality allows them to cover the problems that are non-Lipschitz on $\R^d$ (e.g., for some monotone polynomially-growing operators) and also the situation when the variance is bounded only on a compact, which is common for many finite-sum problems. Finally, \cite{gorbunov2022stochastic} contains high-probability convergence results for strongly monotone VIs and VIs with structured non-monotonicity.

\section{Recent advances}\label{adv}
In this section, we collect a few recent theoretical advances with practical impacts.
\subsection{Saddle point problems with different constants of strong convexity and strong concavity}
Interest in saddle point problems with different constants of strong convexity and strong concavity appeared a few years ago, see e.g., \cite{alkousa2020accelerated,lin2020near}. However, even for the particular case
\begin{equation*}\label{problem_main}
		\min_{x \in \R^{d_x}}\max_{y \in \R^{d_y}}g(x, y) = f(x) + y^\top \mA x - h(y),
\end{equation*}
where function $f$ is $\mu_x$-strongly convex ($\mu_x>0$) and $L_x$-smooth, and function $h$ is $\mu_y$-strongly convex ($\mu_y>0$) and $L_y$-smooth, optimal algorithms have appeared only recently \cite{minimax2112,similar_determ2201,jin2022sharper}.
These algorithms have the following convergence rates
\begin{equation}\label{rate}
\mathcal{O}\left( \left( \sqrt{\frac{L_x}{\mu_x}} + \sqrt{\frac{\lambda_{\max}(\mA^\top\mA)}{\mu_x\mu_y}} + \sqrt{\frac{L_y}{\mu_y}} \right)\log{\frac1\varepsilon} \right),
\end{equation}
and attain the lower bound (here we need to assume that $\lambda_{\min}(\mA^\top \mA) \le \sqrt{\mu_x\mu_y}$, without this assumption optimal methods are  unknown), which was obtained in \cite{zhang2021lower,ibrahim2020linear}

Note that the algorithm from \cite{jin2022sharper} was built upon the technique related to the analysis of
primal-dual \algname{Extragradient} methods via relative Lipschitzness \cite{stonyakin2021inexact,cohen2020relative}. As a by-product, this technique  makes it possible to obtain Nesterov's accelerated method as a particular case of primal-dual \algname{Extragradient} method with relative Lipschitzness \cite{cohen2020relative}.

For the non-bilinear SPP, optimal methods, based on the accelerated Monteiro--Svaiter proximal envelope, were  developed only in the non-composite case \cite{kovalev2022first,carmon2022recapp}. For the non-bilinear SPP with composite terms, there is a poly-logarithmic gap between the lower bound and the best known upper bounds \cite{tominin2021accelerated}. A gap also appears for the SPP having the stochastic finite-sum structure \cite{tominin2021accelerated,luo2021near,jin2022sharper}. The stochastic setting with bounded variance was considered in \cite{couplin_term2111,metelev2022decentralized,du2022optimal}.

 

Further deterministic <<cutting-plane>> improvements are related with the additional assumptions about small dimension of vectors $x$ or/and $y$ \cite{nemirovski1994efficient,gladin2020solving,gladin2021solving} or with different structural (e.g., SPP on balls in $1$/$\infty$-norms) and sparsity assumptions, see e.g.,   \cite{carmon2020coordinate,song2021variance,song2021coordinate} and references there in. Lower bounds here are mostly  unknown. 


Note, that in this subsection we mentioned a lot of works with (sub-)optimal algorithms for different variants of SPP. In contrast to the convex optimization, where the oracle call is uniquely associated with the gradient call $\nabla f(x)$, for SPP we have two criteria: the number of $\nabla_x g(x,y)$-calls and $\nabla_y g(x,y)$-calls (and more variants for SPP with composites). ``Optimality'' in the  most of the aforementioned  papers  means that the method is optimal according to the worst of the criteria. In \cite{alkousa2020accelerated,tominin2021accelerated}, authors consider these criteria separately. However, the development of the lower bounds and optimal methods for a multi-criterion setup is still an open problem.


\subsection{Adaptive methods for VI and SPP}

Interest in adaptive algorithms for  stochastic convex optimization was mainly emerged in 2011 after the development of  \algname{AdaGrad} \cite{duchi2011adaptive} and \algname{Adam} \cite{Kingma2015}. For variational inequalities and saddle point problems, interest in adaptive methods have appeared only in the last few years, see e.g., \cite{gasnikov2019adaptive,bach2019universal} (see also \cite{khobotov1987modification}). Currently, this area of research is well-developed.
One can note works devoted to both adaptive step sizes \cite{antonakopoulos2020adaptive,ene2021adaptive, ene2022adaptive,stonyakin2022generalized,titov2022some} and adaptive scaling/preconditioning \cite{liu2019towards, dou2021one, beznosikov2022scaled}.
Approaches from the second group are based on the idea of the proper combination of  \algname{AdaGrad}/\algname{Adam} with \algname{Extragradient} or its modifications. All of the mentioned adaptive methods have no better (typically the same) theoretical convergence rates than their non-adaptive analogues but require less input information or demonstrate better performance in practice. 

\subsection{Quasi-Newton and tensor methods for VI and SPP}

Quasi-Newton methods for solving nonlinear equations (unconstrained VI) and SPP are proposed in \cite{lin2021explicit, ye2021greedy} and \cite{liu2021quasi} respectively. In these papers, the authors derive local superlinear convergence rates for the modifications of the Broyden's type methods for solving nonlinear equations with Lipschitz Jacobian and SPP with Lipschitz Hessian. Stochastic versions of these methods for VI and SPP have not yet been developed.

Tensor methods for convex optimization problems  are currently quite well developed. In particular, starting with \cite{nesterov2006cubic} it has been  shown that optimal second- and third-order methods can be implemented with almost the same complexity of each iteration as the \algname{Newton method} \cite{monteiro2013accelerated,nesterov2021implementable,gasnikov2019near}. Moreover, optimal $p$-order methods (that use $p$-order derivatives) significantly reduce the rate of convergence from $k^{-2}$ to $ k^{-(3p+1)/2}$ \cite{kovalev2022first_,carmon2022optimal}. For VI and \algname{SPP}, the interest was initiated by \cite{nesterov2006cubic_,monteiro2012iteration} and  optimal $p$-order methods reduce the rate of convergence from $k^{-1}$ to $k^{-(p+1)/2}$ \cite{adil2022line,lin2022perseus} (for $k^{-1}$, see Theorem~\ref{th:MirrProx}). However, in contrast to convex optimization,  the use of tensor methods for sufficiently smooth monotone VIs and convex-concave saddle point problems is not expected to be as effective. Note, that in \cite{adil2022line,lin2022perseus} one can also find optimal rates for strongly monotone VIs and strongly convex-concave saddle point problems. Stochastic tensor methods for variational inequalities and saddle point problems have not yet been developed.

\subsection{Convergence in terms of the gradient norm for SPP}
Several recent advances in  the development of optimal algorithms are based on accelerated proximal  envelops with proper stopping rules for inner loop algorithms \cite{kovalev2022first,kovalev2022first_,kovalev2022optimal_sliding, sadiev2022communication}. This rule is built upon the norm of the gradient calculated for target function of the inner problem.

For smooth convex optimization problems Yu.~Nesterov in 2012 posed the problem of making the gradient norm small with the same rate of convergence as a gap in the function value, i.e. proportional to  $ k^{-2}$ \cite{nesterov2012make}. To address this problem, in the same paper, he proposed an optimal (up to a logarithmic factor) algorithm. This question was further investigated, including obtaining optimal results without additional logarithmic factors \cite{kim2021optimizing,nesterov2021primal} (see also \cite{diakonikolas2022potential} for explanations and survey). In the stochastic case, algorithms were presented in \cite{foster2019complexity}.

For smooth convex-concave saddle point problems an optimal algorithm with $\|\nabla_{x,y} f(x^k,y^k)\|_2$ proportional to $k^{-1}$ was proposed in \cite{yoon2021accelerated} (see also \cite{diakonikolas2022potential} and \cite{kovalev2022first} for monotone inclusion). For the stochastic case see \cite{lee2021fast,cai2022stochastic,chen2022near}.  

\subsection{Decentralized VI and SPP}
 
In practice, in order to solve the variation inequality problem more efficiently and quickly, distributed methods are usually applied. In particular, methods that work on arbitrary (possibly time-varying) decentralized communication networks  between computing devices are popular. 

While the field of decentralized algorithms for minimization problems has been extensively investigated, results for broader classes of problems have only begun to appear in recent years. Such works are primarily focused on saddle point problems \cite{Mukherjee2020:decentralizedminmax, beznosikov2020distributed, rogozin2021decentralized, beznosikov2021distributed, beznosikov2021}, but we note that most of these results can easily be extended to variational inequalities. Let us emphasize two works that were at once devoted to VIs. In the paper \cite{beznosikov2021decentralized} the authors proposed a decentralized method with local steps, and \cite{kovalev2022optimal}  gave optimal decentralized methods for stochastic (finite sum) variational inequalities on fixed and varying networks.

\paragraph{Acknowledgement.}
The work was supported by Russian Science Foundation (project No. 21-71-30005).

\bibliographystyle{plain}
\bibliography{references}

\begin{info}
Aleksandr Beznosikov [anbeznosikov@gmail.com] is a Ph.D. student at Moscow Institute of Physics and Technology (Moscow, Russia). He is also a researcher at Laboratory of Mathematical Methods of Optimization and at Laboratory of Advanced Combinatorics and
Network Applications in Moscow Institute of Physics and Technology (Moscow, Russia), a junior researcher at 
International Laboratory of SA and HDI in Higher School of Economics (Moscow, Russia), a research intern at Yandex Research (Moscow, Russia). His current research interests are concentrated around variational inequalities, saddle point problems, distributed optimization, stochastic optimization, machine learning
and federated learning.

\vspace{0.3cm}

Boris Polyak (1935–2023) was a Head of Ya.Z. Tsypkin Laboratory of Institute for Control Science of Russian Academy of Sciences (Moscow, Russia) and a Professor
of Moscow University of Physics and Engineering (Moscow, Russia). He received a Ph.D. degree in Mathematics from Moscow State University (Moscow, Russia) in 1963 and a Doctor of Science degree (Habilitation) in engineering from Institute for Control Science of Russian Academy of Sciences (Moscow, Russia) in 1977. He authored or coauthored more than 250 papers in peer-review journals and 4 monographs, including "Introduction to Optimization". He was a IFAC Fellow, a recipient of the EURO-2012 Gold Medal and the INFORMS Optimization Society Khyachyan Prize. His main area of research was optimization algorithms and optimal control.

\vspace{0.3cm}

Eduard Gorbunov [ed-gorbunov@yandex.ru] is a researcher Laboratory of Mathematical Methods of Optimization in Moscow Institute of Physics and Technology (Moscow, Russia). His current research interests are concentrated around stochastic optimization and its applications to machine learning, distributed optimization, derivative-free optimization, and variational inequalities.

\vspace{0.3cm}

Dmitry Kovalev [dakovalev1@gmail.com] is a Ph.D. student at King Abdullah University of Science and Technology (Thuwal, Saudi Arabia). In 2020 and 2021 he received the Yandex Award (Ilya Segalovich Award). His current research interests include continuous optimization and machine learning.

\vspace{0.3cm}

Alexander Gasnikov [gasnikov@yandex.ru] is a professor at Moscow Institute of Physics and Technology, head of the Laboratory Mathematical methods of optimization and head of the department Mathematical foundations of control. He received a Doctor of Science degree (Habilitation) in Mathematics in 2016 from the Faculty of Control and Applied Mathematics of Moscow Institute of Physics and Technology. In 2019 he received Yahoo Faculty Research Engagement Program. In 2020 he received the Yandex Award (Ilya Segalovich Award). In 2021 he received the Award for Young Scientists from the Moscow Government. His main area of research is optimization algorithms.

\end{info}

\end{document}